\documentclass[12pt,twoside,leqno,letterpaper]{article}
\usepackage{amsmath,amssymb,multicol}
\newcommand{\Finalday}{23rd October 2006} \date{\Finalday}
\makeatletter
\newcommand{\provsection}{\@startsection%
{section}
{1}
{\parindent}
{-\baselineskip}
{-0.5\baselineskip}
{\normalfont\Large\bfseries}
}%
\makeatother%
\renewcommand{\section}[1]{\provsection{#1}}
\makeatletter
\newcommand{\provsubsection}{\@startsection%
{subsection}
{2}
{\parindent}
{-\baselineskip}
{-0.5\baselineskip}
{\normalfont\normalsize\bfseries}
}%
\makeatother%
\renewcommand{\subsection}[1]{\provsubsection{#1.}}
\makeatletter
\newcommand{\provsubsubsection}{\@startsection%
{subsubsection}
{3}
{\parindent}
{-0.2\baselineskip}
{-0.5\baselineskip}
{\normalfont\normalsize\bfseries}
}%
\makeatother%
\renewcommand{\subsubsection}[1]{\provsubsubsection{#1}}
 
\renewcommand{\[}{\begin{eqnarray*}}
\renewcommand{\]}{\end{eqnarray*}}
\newcommand{\la}{\begin{eqnarray}}
\newcommand{\al}{\end{eqnarray}}

\newcommand{\conv}{\ast}

\DeclareMathOperator*{\bigconv}{\mbox{\LARGE$\conv$}}
\newcommand{\Halmos}{\mbox{}\hfill\rule{2mm}{2mm}}

\renewcommand{\d}{{\,\text{\rm d}}}
\newcommand{\e}{\mathrm{e}}
\newcommand{\pii}{\mathrm{\pi}}

\newcommand{\BesselI}{\mathrm{I}}                        %
\newcommand{\binomial}[2]{\genfrac{(}{)}{0pt}{}{#1}{#2}} 

\newcommand{\set}[1]{{\left\{#1\right\}}}

\renewcommand{\epsilon}{\varepsilon}\renewcommand{\phi}{\varphi}      %
\renewcommand{\rho}{\varrho}        \renewcommand{\theta}{\vartheta}  %
\newcommand{\C}{{\mathbb C}}

\newcommand{\N}{{\mathbb N}}
\newcommand{\R}{{\mathbb R}}
\newcommand{\Z}{{\mathbb Z}}

\newcommand{\cP}{{\cal P}}%
\newcommand{\1}{\mathbf{1}}
\renewcommand{\P}{{\mathbb P}}                                               %
\newcommand{\Conc}{\mathrm{Conc}}                                            %
\newcommand{\Ber}{\mathrm{B}}     
\newcommand{\BerC}{\mathrm{B}}    
\newcommand{\Bin}{\mathrm{B}}     
\newcommand{\Pois}{\mathrm{P}}    
\newcommand{\STP}{\mathrm{Q}}     
\newcommand{\STPC}{\mathrm{Q}}    
\newcommand{\RadC}{\mathrm{R}}    
\newcommand{\SymPois}{\mathrm{S}} 

\newcommand{\pollard}{ \sim\kern -0.7ex>}

\title{A shorter proof of \\
 Kanter's Bessel function concentration bound
        }

\author{Lutz Mattner \& Bero Roos \\
         Universit\"at zu L\"ubeck  \& Universit\"at Hamburg  }

\begin{document}
\maketitle
\begin{abstract}
We give a shorter proof of Kanter's (1976) sharp Bessel function bound 
for concentrations of sums 
of independent symmetric random vectors. 
We provide sharp upper bounds for the sum
of modified Bessel functions 
$\BesselI_0(x)+\BesselI_1(x)$, which might
be of independent interest.
Corollaries improve  concentration or smoothness  bounds 
for sums of independent random variables due
to \v{C}ekanavi\v{c}ius \& Roos~(2006), Roos~(2005), 
Barbour \& Xia~(1999), and Le~Cam~(1986).   
\end{abstract}



\section{Introduction} \mbox{} \\

The principal purpose of this paper is to shorten one proof
in the development of Kanter's (1976) remarkable concentration 
bounds for sums of independent random vectors. Restricting our 
attention to the symmetric and finite-dimensional case,
Kanter's  main results may be stated as follows.

Let us call a random vector $X$ symmetrically distributed
if $X$ and $-X$ have the same distribution.
Let  $\|\cdot\|$ be a seminorm on a finite-dimensional 
$\R$-vector space $E$,  let $X_1,\ldots,X_n$ 
be independent and symmetrically distributed $E$-valued random vectors,
let $t\in\,]0,\infty[$, and let us put
\la \label{Def pj from Xj} 
  p_j &:=& \P(\|X_j\| \ge t) \qquad (j=1,\ldots,n)
\al 
Let $Y_1,\ldots,Y_n$ be independent and symmetrically
distributed $\R$-valued random variables with 
\la                 \label{Distr of Yj}
  \P(|Y_j|=1)  &=&1-  \P(Y_j=0)  \,=\, p_j \qquad (j=1,\ldots,n)
\al
Let further
\la                   \label{Def G neu}
 G(\lambda) &:=& 
  \mathrm{e}^{-\lambda}\big(\BesselI_0(\lambda) +\BesselI_1(\lambda)\big)
 \qquad (\lambda \in [0,\infty[)
\al
where $\BesselI_k$ denotes the modified Bessel function of order $k$. 
Then  the following three results hold.
\subsection{Theorem (Kanter's Lemma 4.2)}\label{Theorem (Kanter's Lemma 4.2)}
{\sl
\la    \label{Ineq Theorem (Kanter's Lemma 4.2)}
  \P(\| \sum_{j=1}^n X_j - x \| < t) &\le& 
  \P(| \sum_{j=1}^n Y_j - \frac 12  | < 1)  \qquad (x\in E)
\al}%
  
\subsection{Theorem (Kanter's Lemma 4.3)}\label{Theorem (Kanter's Lemma 4.3)}
{\sl If\, $p_j>0$ for some $j$, then  
\la \label{Ineq Theorem (Kanter's Lemma 4.3)}
 \P(| \sum_{j=1}^n Y_j - \frac 12  | < 1) 
  &<& G(\sum_{j=1}^n p_j)
\al}%

\subsection{Corollary (Kanter's Theorem 4.1)}
\label{Corollary (Kanter's Theorem 4.1)}
{\sl If\,  $\P(|X_j\| <t) <1$ for some $j$, then
\la  \label{Ineq Corollary (Kanter's Theorem 4.1)}
\sup_{x\in E}  \P(\| \sum_{j=1}^n X_j - x \| < t) &<& 
   G\Big(\sum_{j=1}^n\big(1-\P(\|X_j\| < t)\big)\Big)
\al}

All three results are optimal: 
In Theorem \ref{Theorem (Kanter's Lemma 4.2)}, equality obviously holds
whenever $x_0\in E$ has norm $t$, $x=x_0/2$, and 
$\P(X_j=x_0)=\P(X_j=-x_0) = p_j/2 = \big( 1- \P(X_j=0)\big)/2 $
for each $j$. Thus Theorem \ref{Theorem (Kanter's Lemma 4.2)}
provides the supremum  of the left hand side
of \eqref{Ineq Theorem (Kanter's Lemma 4.2)} given 
the probabilities $p_j$ from \eqref{Def pj from Xj}.
See Mattner (2006) for a generalization.
In Theorem \ref{Theorem (Kanter's Lemma 4.3)},
the left hand side of \eqref{Ineq Theorem (Kanter's Lemma 4.3)},
being some  rather complicated function of 
$(p_1,\ldots,p_n) \in \bigcup_{n\in\N} [0,1]^n$,
is bounded by its supremum given $\sum_1^n p_j$, 
see Remark  \ref{Remarks in its own right} \ref{Sharpness}
below. It follows that 
Corollary \ref{Corollary (Kanter's Theorem 4.1)}
provides the supremum of the left hand side in 
\eqref{Ineq Corollary (Kanter's Theorem 4.1)} given 
$\sum_{j=1}^n\big(1-\P(\|X_j\| < t)\big)$. 

To further simplify the bounds from  
\eqref{Ineq Theorem (Kanter's Lemma 4.3)}
or \eqref{Ineq Corollary (Kanter's Theorem 4.1)},
one may use the inequalities from \eqref{simple bound, 1},
\eqref{simple bound, 2} or \eqref{simple bound, 3}
of the following collection of analytical properties of the function $G$.

\subsection{Lemma}\label{Lemma Simple bounds for G}
{\sl $G$ is completely monotone on $[0,\infty[$,
with $G(0)=1$ and $\displaystyle\lim_{\lambda\rightarrow \infty}G(\lambda)=0$. 
For $\lambda \in\,]0,\infty[$, we have
\la
  G(\lambda) &=&  \sqrt{\frac 2{\pi \lambda}}  
   \Big(1 -\frac 1{8\lambda} -\frac 3{128 \lambda^2}  
  + O(\frac 1{\lambda^3}) \Big) \label{G asymp for lambda large} \\  
  G(\lambda) &<& \sqrt{\frac 2{\pi \lambda}} \label{simple bound, 1} \\
G(\lambda) &<& \sqrt{\frac {\frac 2\pi}{\frac 14 +\lambda}}
    \,=\, \sqrt{\frac 2{\pi\lambda}}\Big(1-\frac 1{8\lambda} 
       + \frac 3{128\lambda^2} 
 +O(\frac 1{\lambda^3})\Big) \label{simple bound, 2} \\
  G(\lambda) &=& 1-\frac{\lambda}2 +\frac{\lambda^2 }4+O(\lambda^{3})
                       \label{G asymp for lambda small} \\
G(\lambda) &<& \frac 1{\sqrt{1+\lambda}}  
\,=\,
1-\frac{\lambda}2 +\frac{3\lambda^2}{8}+O(\lambda^3) \label{simple bound, 3}
\al}

The present paper is mainly  concerned with 
Theorem \ref{Theorem (Kanter's Lemma 4.3)}. We found 
Kanter's original proof rather long and somewhat 
difficult to read, while  a much shorter proof of a stronger 
claim proposed by Marshall \& Olkin (1979) is unfortunately not
generally valid, see 
Remark  \ref{Remarks in its own right}\ref{Remark on Marshall and Olkin}
below. Our proof of Theorem \ref{Theorem (Kanter's Lemma 4.3)}
given below
is shorter than Kanter's,  and we hope that it is
reasonably  easy to read.  As Theorem 
\ref{Theorem (Kanter's Lemma 4.3)}  can be read
as an inequality between certain special convolution products 
and their limits,
we have chosen to give an introduction from that point 
of view  in 
Section~\ref{Symmetric three point convolutions and related distributions}
below, where Theorem~\ref{Theorem (Kanter's Lemma 4.3)} reappears 
as the slightly more general
Theorem~\ref{Kanter Lemma 4.3}. To check that
the latter implies the former,   use
\eqref{Kanter Lemma 4.3 inequality} with $k=0$, 
\eqref{Glambda vs Slambda01}, and the fact  that $\sum_{j=1}^n Y_j$ from 
\eqref{Ineq Theorem (Kanter's Lemma 4.3)} has the distribution
$\STPC_p$ introduced at the beginning of 
Section~\ref{Symmetric three point convolutions and related distributions}.
This latter fact together with  \eqref{ConcQp as Bernoulli conv expectation}
also shows that Theorem 
\ref{Theorem (Kanter's Lemma 4.2)}
is indeed a reformulation of Kanter's Lemma 4.2.

Theorem~\ref{Kanter Lemma 4.3} is proved 
in  \ref{Proof of  Kanter's Lemma 4.3} below,
relying on  a crucial  analytic inequality
provided in \ref{special analytic inequality}.
The rest of Section \ref{Proofs} is then devoted to proofs of 
Lemma~\ref{Lemma Simple bounds for G} 
and  Corollaries~\ref{Corollary improving Roos 2005} 
and \ref{Corollary improving Barbour Xia 1999} stated below,
which are not needed for the proof of Theorem~\ref{Kanter Lemma 4.3}.

It is not  the purpose of this paper to treat 
systematically concentration function bounds  related to Kanter's.
See Petrov (1995) for a good introduction to the one-dimensional case. 
A complete review  would include in particular
Bretagnolle (2004)  and  many references from the journal
{\it Theory  of Probability and Its Applications}, such as Rogozin (1993).
Let us however give  the following corollary to  
Theorem \ref{Theorem (Kanter's Lemma 4.3)}
and  Lemma \ref{Lemma Simple bounds for G}
on concentration functions $\Conc(X,\cdot)$ 
of real-valued random variables $X$,
defined by 
\[
\Conc(X,t) &:=& \sup_{x\in\R} \P( X \in [x,x+t]) \qquad (t\in[0,\infty[)
\]
We note that \eqref{bound with sharp constant} below
complements  \v{C}ekanavi\v{c}ius \& Roos (2006, Lemma 3.6),
while \eqref{bound improving Roos 2005}
simplifies and uniformly improves Roos (2005, Proposition 5) 
and hence as well Le~Cam (1986, p.~411, Theorem 2), 
and complements the results of Bretagnolle. 
Let $\1_A(x)$ be $1$ or $0$
according to whether $x\in A$ or $x\notin A$.
Random variables $X_j$ are no longer supposed to be symmetric 
and numbers $p_j$ may differ from those defined in \eqref{Def pj from Xj}. 

\newcommand{\otherp}{p}
\subsection{Corollary} \label{Corollary improving Roos 2005}
{\sl Let $X_1,\ldots,X_n$ be independent 
$\R$-valued random variables and let $t\in[0,\infty[$. 
For $j=1,\dots,n$, let
$h_j(y):=\inf\{x\in\R\,:\,\P(X_j\leq x)\geq y\}$ for  $y\in\,]0,1[$
and $\otherp_j:=2\int_0^{1/2}\1_{]t,\infty[}(h_j(1-y)-h_j(y))\d y$.
Then  
$\otherp_j\geq 1-\Conc(X_j,t)$  for every $j$ and
\la
 \Conc (\sum_{j=1}^n X_j, t)  
 &\le & 
 G\big(\sum_{j=1}^n \otherp_j \big) 
 \,\, \le\,\,
 G\Big(\sum_{j=1}^n \big(1-\Conc(X_j,t)\big) \Big) 
\al
and hence in particular
\la
\Conc (\sum_{j=1}^n X_j, t) 
 &\le& \sqrt{\frac 2\pi}\Big(\frac 14 +\sum_{j=1}^n \otherp_j  \Big)^{-1/2}
 \label{bound with sharp constant}\\
 \Conc (\sum_{j=1}^n X_j, t)  
 &\le& \Big( 1+\sum_{j=1}^n \big(1-\Conc(X_j,t)\big)\Big)^{-1/2}
 \label{bound improving Roos 2005}
\al}%

In \eqref{bound with sharp constant}, 
the constant $\sqrt{\frac{2}{\pii}}$ is optimal, since for
$\P(X_j=-1)=\P(X_j=1)=1/2$, we have $\otherp_j=1$
and $\Conc(\sum_{j=1}^n X_j, 1)  
\sim \sqrt{2/\pii}(\sum_{j=1}^n \otherp_j)^{-1/2}$ for $n\rightarrow\infty$. 

Theorem~\ref{Theorem (Kanter's Lemma 4.3)}
and Lemma~\ref{Lemma Simple bounds for G} 
further yield an improvement of a
smoothness bound of Barbour \& Xia (1999, Proposition~4.6),
who applied it in the context of compound Poisson approximation.
We denote the total variation distance between distributions
of $\Z$-valued random variables $Z_1$ and~$Z_2$ by 
$d_{\mathrm{TV}}(Z_1,\,Z_2):=\sup_{A\subseteq\Z}|\P(Z_1\in A)-\P(Z_2\in A)|$.
\subsection{Corollary} \label{Corollary improving Barbour Xia 1999}
{\sl Let $X_1,\ldots,X_n$ be independent $\Z$-valued random variables.
Then
\[
 d_{\mathrm{TV}}  (\sum_{j=1}^n X_j,\, 1+\sum_{j=1}^n X_j)
 &\le& \sqrt{\frac 2\pi}
 \Big(\frac 14 +\sum_{j=1}^n\big(1-d_{\mathrm{TV}}(X_j,\,1+X_j)\big)
 \Big)^{-1/2}
\]} 

This improves Barbour \& Xia's bound 
$ (\sum_{j=1}^n
(1-\max(\frac 12, d_{\mathrm{TV}}(X_j,\,1+X_j))))^{-1/2}$.

\section{Symmetric three point convolutions and related distributions}
\label{Symmetric three point convolutions and related distributions}
\mbox{}\\
 
We let $\N:= \set{1,2,3,\ldots}$ and $\N_0:= \set{0} \cup\N$.
Throughout the paper, we assume $n\in\N_0$ and 
$p=(p_1,\ldots,p_n)\in \cP:= \bigcup_{n\in\N_0}[0,1]^n$
and put
$|p| := \sum_{j=1}^n p_j$,  
$\overline{p} := |p|/n$, $p_{\max} := \max_{j=1}^n p_j$,
and $n(p):=n$. 
For $\alpha\in[0,1]$, 
$p
\in \cP $ 
and $\lambda\in[0,\infty[$,
we consider the following probability  distributions on $\Z$: 
\newcommand{\phant}{  \displaystyle\phantom{\sum_1}}
\[
\begin{array}{lrcl} 
 \mbox{Bernoulli:} & \Ber_\alpha &:=& 
     (1-\alpha)\delta_0+\alpha\delta_1 \phant\\
 \mbox{Bernoulli convolution:}& \displaystyle
 \BerC_{p} &:=&  \displaystyle\bigconv_{j=1}^n \Ber_{p_j}  \phant \\
 \mbox{Binomial:}&\displaystyle 
 \Bin_{n,\alpha} &:=&  \Ber_\alpha^{\conv n} \,=\, 
       \BerC_{(\alpha,\ldots,\alpha)} \phant   \\
 \mbox{Poisson:}&\displaystyle
  \Pois_\lambda &:=& 
 \displaystyle\sum_{j=0}^\infty \e^{-\lambda}\frac{\lambda^j}{j!}\delta_j  \\
 \mbox{Symmetric $3$-point:}
 &\displaystyle
 \STP_\alpha &:=& 
  (1-\alpha)\delta_0
   +\frac\alpha2\big(\delta_{-1}+ \delta_1\big) \phant\\
 \mbox{Symmetric $3$-pt.\ convolution:} &\displaystyle
 \STPC_p &:=& \displaystyle\bigconv_{j=1}^n \STP_{p_j} \phant\\
 \mbox{Rademacher convolution:}  &\RadC_n &:=& \displaystyle
    \Big(\frac{\delta_{-1}+\delta_1}2 \Big)^{\conv n}
  \,=\, \STPC_{(1,\ldots,1)}
   \,=\,  \sum_{j=0}^n \binomial{n}{j}  2^{-n}\delta^{}_{2j-n} \phant \\
 \mbox{Symmetrized Poisson:} 
 &\SymPois_\lambda &:=& \Pois_{\lambda/2}\conv \widetilde{\Pois}_{\lambda/2}
 \,=\,\displaystyle \sum_{j\in \Z} \e^{-\lambda}\BesselI_j(\lambda)\delta_j
        \phant  \\  
\end{array}
\]
We use roman letters for these special probability distributions,
reserving the italic letters $P$ and $Q$  for arbitrary distributions  
on $\Z$.  Here and in what follows, $\widetilde{P} $ denotes the reflection of
 $P$, 
defined by $\widetilde{P}(\set{k}):= P(\set{-k})$, $\delta_k$
is the Dirac measure concentrated in $k$, $\conv $ denotes convolution,
and empty convolution products and powers are understood to be  $\delta_0$. 
The last representation for the symmetrized Poisson distributions 
$\SymPois_\lambda$ follows from  the power 
series expansion of the modified Bessel functions $\BesselI_j$, namely 
$\BesselI_j(x)=\BesselI_{|j|}(x) 
= \sum_{k=0}^\infty \frac 1{k!(k+|j|)!}(\frac x2)^{2k+|j|}$ for 
$j\in\Z$ and $x\in \C$, see for example 
Olver (1997, p.~60). Recalling \eqref{Def G neu}, we observe that 
\la \label{Glambda vs Slambda01}
   G(\lambda) &=& \SymPois_\lambda(\set{0,1}) \qquad (\lambda \in[0,\infty[)
\al

The  objects of our study are the symmetric $3$-point convolutions 
$\STPC_p$. Let us collect some simple properties of these. 
For every fixed $n$, the probability measure   $\STPC_p$ is a 
permutation invariant and multiaffine function
of the parameter $p\in[0,1]^n$. Hence, 
as noted by Hoeffding (1956, p.~713),  each $\STPC_p$ admits a 
 representation as an  expectation with respect to a unique
Bernoulli convolution $\BerC_p$. In fact, if we put
$P:= \delta_0$ and $Q:= \frac 12(\delta_{-1}+\delta_1)$
in the identity 
\[ 
\bigconv_{j=1}^n\big((1-p_j)P+p_jQ\big)
 &=& \sum_{m=0}^n \Ber_p(\set{m})Q^{\conv m} \conv P^{\conv(n-m)}
\]
we get
\la \label{Q as Rademacher conv mixture} 
 \STPC_p &=& \int_{\N_0} \RadC_m \d \BerC_p(m) \qquad(p\in\cP)
\al
as a mixture of Rademacher convolutions.
In case of $p_{\max} \le1/2$,
the symmetric $3$-point convolution $\STPC_p$   
is also a symmetrized Bernoulli convolution, since 
\[
  \STP_\alpha &=& 
 \Ber^{}_{\beta(\alpha)} \conv \widetilde{\Ber}^{}_{\beta(\alpha)} 
 \quad\text{ for }\quad\alpha \in[0,1/2]\quad\text{ and }\quad  
 \beta(\alpha) := \frac {1-\sqrt{1-2\alpha}}{2}  
\]
and hence 
\la
 \STPC_p &=& \Big( \bigconv_{j=1}^n \Ber^{}_{\beta(p_j)}\Big) \conv  
  \Big(\bigconv_{j=1}^n \Ber^{}_{\beta(p_j)}\Big)^{\widetilde{\rule{5pt}{0pt}}}
 \qquad (n\in\N_0,\,p\in [0,1/2]^n)
\al
Starting  from this representation, one can easily show
for every  $\lambda \in [0,\infty[$ that $\STPC_p$ converges 
to $\SymPois_\lambda$   in  total variation 
if $p\in \cP $  
varies in such a way that $|p| \rightarrow \lambda$
and $p_{\max} \rightarrow 0$:
\la       \label{lim Qp=Slambda}
  \lim_{|p| \rightarrow\lambda,\,
       p_{\max} \rightarrow 0} \STPC_p &=& \SymPois_\lambda
  \qquad (\lambda \in[0,\infty[)
\al 
To verify this, one may use  a classical Poisson
approximation theorem for Bernoulli convolutions, 
for example 
Barbour {\it et al.} (1992, p.~3,  Le Cam's result  (1.6)),
the fact that  $\beta(p_j) =\big(1+ O(p_{\max})\big)\frac{p_j}2$ 
uniformly in $j$, and the continuity of the map 
$\lambda \mapsto  \Pois_\lambda$.
Below we will also need the representation
\la \label{S as Rademacher conv mixture} 
 \SymPois_\lambda 
 &=& \int_{\N_0} \RadC_m \d \Pois_\lambda(m) \qquad(\lambda\in[0,\infty[)
\al
This  can be proved similarly to \eqref{Q as Rademacher conv mixture}   
by  a short computation or, alternatively, by combining  
 \eqref{Q as Rademacher conv mixture} with \eqref{lim Qp=Slambda}.

Theorem \ref{Theorem (Kanter's Lemma 4.3)}
is related to, but  decidedly  less obvious   
than \eqref{lim Qp=Slambda}. Rephrased with the notation of this section,
it states that   $\SymPois_\lambda$ with $\lambda = |p|$ even 
serves as a bound for, and not merely as an approximation to
$\STPC_p$, as far as the concentration over two-point
intervals in $\Z$  is concerned.

\subsection{Theorem (= Theorem \ref{Theorem (Kanter's Lemma 4.3)})} 
\label{Kanter Lemma 4.3}
{\sl Let $p\in \cP$ with $|p|>0$ and let $k\in \Z$.
Then
\la   \label{Kanter Lemma 4.3 inequality}
 \STPC_p(\set{k,k+1}) &<& \SymPois_{|p|}(\set{0,1})   
\al}%

\subsection{Remarks} \label{Remarks in its own right}
\subsubsection{} \label{Sharpness}
Inequality \eqref{Kanter Lemma 4.3 inequality}
is sharp: For every fixed $\lambda >0$,
the supremum of the left hand side of \eqref{Kanter Lemma 4.3 inequality}
over all $p\in\cP$ with $|p|=\lambda$ and all $k$ is 
$\SymPois_\lambda(\set{0,1})$, as  follows using \eqref{lim Qp=Slambda}.

\subsubsection{} \label{Where max?}
As a function of $k\in\Z$, each of 
$\STPC_p(\set{k,k+1})$  and $\SymPois_\lambda(\set{k,k+1})$
becomes maximal iff $k=0$ or $k=-1$.
This follows from the mixture representations 
\eqref{Q as Rademacher conv mixture} and \eqref{S as Rademacher conv mixture},
since $\RadC_m(\set{k,k+1})$ obviously becomes maximal iff $k=0$ or $k=-1$.
Hence we may rewrite Theorem \ref{Kanter Lemma 4.3} as the 
concentration function inequality
\la                \label{Kanter special conv.ineq}
 \sup_{k\in\Z} \STPC_p(\set{k,k+1}) 
  &<& \sup_{k\in\Z} \SymPois_{|p|}(\set{k,k+1}) 
  \qquad(p\in \cP,\,|p| >0)
\al 
\subsubsection{} \label{unimodality}
The distributions $\SymPois_\lambda$ are  symmetric and, if $\lambda >0$, 
strictly unimodal on $\Z$, that is, we have
$\SymPois_\lambda(\set{k})=\SymPois_\lambda(\set{|k|})
>\SymPois_\lambda(\set{|k|+1})$
for every $k\in \Z$.  Here  the claimed inequality
is a special case of the known strict antitonicity 
in the order $\nu \in\,]0,\infty[$ of the 
modified Bessel functions $\BesselI_\nu(x)$ at  fixed 
 arguments $x\in\,]0,\infty[$,
see  Olver (1997, p.~251, Theorem 8.1(ii)).
For   $p \in [0,2/3]^n$,  each $\STP_{p_j}$ is
symmetric and unimodal on $\Z$, and hence so is their
convolution $\STPC_p$, 
by the discrete Wintner theorem in
Dharmadhikari \& Joag-Dev~(1988, p.~109, Theorem 4.7).

\subsubsection{} \label{ell-Punkt statt 2-Punkt}
Inequality  \eqref{Kanter special conv.ineq}
for two-point intervals in $\Z$ does not generalize in the
obvious way to more general $\ell$-point intervals $\set{k,\ldots,k+\ell -1}$
with $\ell \in\N$. More precisely, there is  no 
$\ell \in\N\setminus\set{2}$ such that  the inequality
\la       \label{Kanter special conv.ineq generalized}
 \sup_{k\in\Z} \STPC_p(\set{k,\ldots,k+\ell-1}) 
  & \le & \sup_{k\in\Z} \SymPois_{|p|}(\set{k,\ldots,k+\ell-1}) 
\al
holds for every $p\in\cP$. 
For $\ell = 1$, we get a counterexample to 
\eqref{Kanter special conv.ineq generalized}
by  choosing $p=(1,\ldots,1)\in\set{1}^n$
with  $n:=2m$, for large enough $m\in\N$.
In fact, for $m\rightarrow \infty$, we have 
$\STPC_p(\set{0})=\RadC_{2m}(\set{0}) =\binom{2m}{m}2^{-2m} 
\sim 2/\sqrt{2\pi n}$ 
but, using Remark \ref{unimodality}, we have
$\sup_{k\in\Z}\SymPois_n(\set{k})= \SymPois_n(\set{0})
=\e^{-n}\BesselI_0(n) \sim 1/\sqrt{2\pi n}$.
See  Olver (1997, p.~83) or \eqref{BesselI for large argument}         
below for the standard Bessel function asymptotics just used.
For $\ell \ge 3$, we get a trivial counterexample 
by choosing $p \in \,]0,1]$ one-dimensional,
since then the left hand side of 
\eqref{Kanter special conv.ineq generalized}
equals  $1$.  More generally,
\[
G^{}_\ell(\lambda) &:=&  
\sup \Big\{\sup_{k\in\Z} \STPC_p(\set{k,\ldots,k+\ell-1}) \,:\,
     p\in \cP, |p| =\lambda      \Big\} 
     \qquad   (\ell\in\N,\,\lambda \in[0,\infty[)
\]
equals $1$ for $\lambda \le \lfloor (\ell-1)/2 \rfloor $,
where, as usual, $\lfloor x \rfloor := \sup \set{n\in\Z : n\le x}$.
While Theorem~\ref{Kanter Lemma 4.3} and Remark~\ref{Sharpness}
state that $G_2$ is just $G$ from \eqref{Def G neu} and 
\eqref{Glambda vs Slambda01}, 
it is an open problem to  compute $G_\ell(\lambda)$
for  $\ell \in\N\setminus\set{2}$ and 
$\lambda >\lfloor (\ell-1)/2 \rfloor$.

\subsubsection{}
Applying  \eqref{Q as Rademacher conv mixture} and 
\eqref{S as Rademacher conv mixture} to the event $\set{0,1}$, and writing 
\la \label{Def psi}
  \psi(m) &:=& \RadC_m(\set{0,1}) 
  \,=\, {m\choose {\lfloor \frac {m+1}2  \rfloor}}2^{-m}
  \,=\, \binomial{2\lfloor \frac {m+1}2\rfloor}{\lfloor \frac {m+1}2\rfloor}
         2^{-2\lfloor \frac {m+1}2\rfloor}
 \qquad (m\in\N_0) 
\al
we get, using the notation $P\phi := \int \phi \d P$,  
\la     \label{ConcQp as Bernoulli conv expectation}
  \STP_p(\set{0,1}) 
   &=& 
  \BerC_p \psi \qquad (p\in\cP)  \\
  \SymPois_\lambda(\set{0,1})  &=& \Pois_\lambda \psi
       \qquad (\lambda\in[0,\infty[)
\al   
so that Theorem~\ref{Kanter Lemma 4.3} can also be stated as 
\la                 \label{Kanter Lemma 4.3 as BerC bound}
   \BerC_p \psi &<& \Pois_{|p|} \psi  \qquad (p\in\cP,\, |p|>0)
\al
Using most conveniently the last expression  for $\psi$ from 
\eqref{Def psi}, we compute
\la
                           \label{special values psi for 0,1}
\qquad \psi(0)=1, \qquad \psi(1)=\psi(2) = \frac12, 
            \qquad \psi(3)=\psi(4) = \frac 38
\al

\subsubsection{} \label{Remark on Marshall and Olkin}
Marshall \& Olkin (1979, p.~378, Theorem K.4, the case $m=1$)
claim   that $\STP_p(\set{0,1})$ 
is, for fixed $n\in \N$,  a Schur-concave function of $p\in[0,1]^n$.
Their proof and claim  become correct if ``$p\in[0,1]^n$'' 
is replaced by  ``$p\in[0,1/2]^n$''.
Without this change  their  claim is  
false.  To see this,  let us note that the following 
three properties of
functions $\phi :\N_0\rightarrow \R$ are equivalent: 

(i) $\phi$  is convex, that is, 
$\phi(k+2)-2\phi(k+1)+\phi(k) \ge 0$ holds for 
$k\in \N_0$, 

(ii) for every $n\in \N$, the Bernoulli convolution expectation 
$\BerC^{}_p \phi$ is a Schur concave function
of $p\in[0,1]^n$,

(iii) $ \BerC^{}_p \phi \le \Bin^{}_{n,\overline{p}}\phi$
for every $n\in \N$ and $p\in[0,1]^n$.

Here the implication ``(i) $\Rightarrow$ (ii)'' 
is due to Karlin \& Novikoff (1963, pp.~1257--1258),
``(ii) $\Rightarrow$ (iii)'' is clear,
and ``(iii) $\Rightarrow$ (i)'' can be seen as follows:
For given $k\in \N$, let  
$n:= k+2$, $\epsilon \in \,\,]0,1/2]$,  and
$p:=(1-\epsilon,\ldots,1-\epsilon,1,1-2\epsilon) \in[0,1]^n$,
so that   
$\BerC^{}_p = 
 \Bin^{}_{k,1-\epsilon}\conv\delta^{}_1\conv\Ber^{}_{1-2\epsilon}$.
Then  $\overline{p}=1-\epsilon$ and 
$\Bin^{}_{n,\overline{p}}- \BerC^{}_p  
= \epsilon ^2 \Bin^{}_{k,1-\epsilon} \conv (\delta_0-2\delta_1+\delta_2) 
 = \epsilon^2(\delta_{k}-2\delta_{k+1} + \delta_{k+2}) +O(\epsilon^3) $,
so that  $\BerC^{}_p \phi \le \Bin^{}_{n,\overline{p}} \phi$
for every $\epsilon$ implies $\phi(k+2)-2\phi(k+1)+\phi(k) \ge 0$.
We remark  that,  similarly,   (ii) with $n$ fixed
implies the convexity of $\phi|_\set{0,\ldots,n}$.

We refer to Hoeffding (1956, Theorem 3) for the original 
proof of ``(i) $\Rightarrow$ (iii)'', to Gleser (1975) 
and Boland \& Proschan (1983) for
generalizations of ``(i) $\Rightarrow$ (ii)'' under restrictions
on $p$, and to Bickel \& van Zwet (1980) and 
Berg {\it et al.} (1984, Chapter 7) 
for results involving distributions more general 
than Bernoulli convolutions.

Coming back to the claim of Marshall \& Olkin,
we observe that the function $\psi$ from \eqref{Def psi}  and 
\eqref{special values psi for 0,1} is not convex on $\set{0,1,2,3}$,
so that it follows, via \eqref{ConcQp as Bernoulli conv expectation}
and using the remark after the above proof of  ``(iii) $\Rightarrow$ (i)'',
that  $\STP_p(\set{0,1})$ is not a Schur concave function 
of $p\in[0,1]^n$ in case $n\ge 3$.

Let us finally note  that Bondar (1994) gives valuable  comments 
on and corrections to Marshall \& Olkin (1979) in general, 
but does not refer to their claim cited above,
while Merkle \& Petrovi\'c (1997, p.~175) repeat the 
claim and attribute it to Kanter (1976).

\subsubsection{} The Poisson bound $\Pois^{}_{|p|} \psi$ on the right 
in \eqref{Kanter Lemma 4.3 as BerC bound} can not in general
be replaced by the  binomial bound $\Bin^{}_{n(p),\overline{p}}\psi$. 
This follows from the implication  ``(iii) $\Rightarrow$ (i)'' 
of the preceding remark.

\section{Proofs} \label{Proofs}
\subsection{A special analytic inequality}\label{special analytic inequality}
{\sl Let $\lambda\in \,]0,\infty[$,  $\alpha \in \,]0,1]$, and
\la
 F(\lambda,\alpha) &:=&      \label{Def F alpha} 
  \frac 1\pi \int_0^\pi 
  \big|1-\alpha(1-\cos t)\big|^{\frac \lambda \alpha} 
     (1+\cos t)\d t 
\al
Then
\la                         \label{crucial ineq}
  F(\lambda,\alpha)&<& G(\lambda)
\al}%

{\bf Remark.} 
It easily follows  that
$\sup_{\alpha\in \,]0,1]} F(\lambda,\alpha) = 
\lim_{\alpha\downarrow 0} F(\lambda,\alpha) = G(\lambda)$,
using \eqref{G by Fourier}.

{\bf Proof.} By Olver (1997, p.~60) or by \eqref{Glambda vs Slambda01} 
and  Fourier inversion of $\SymPois_\lambda$, we have
\la                            \label{G by Fourier}
 G(\lambda)  &=&     \label{Def G}  
  \frac 1\pi \int_0^\pi \e^{-\lambda(1-\cos t)}(1+\cos t)\d t
\al

If $\alpha \in \,\,] 0,1/2]$, then 
$  \big|1-\alpha(1-\cos t)\big|
 = 1-\alpha(1-\cos t) 
 < \exp \big(-\alpha (1-\cos t)\big)$
for every $t\in \,]0,\pi]$, yielding \eqref{crucial ineq}.

So let $\alpha \in \,\,]\frac 12,1]$. 
We split the integral from \eqref{Def F alpha} 
as $\int_0^{t(\alpha)} + \int_{t(\alpha)}^\pi$ with  
$t(\alpha) := \arccos (1- \frac 1\alpha)$,  substitute 
$t = \arccos \big( 1-\frac 1\alpha(1-\e^{- x})\big)$ in the first
integral, and    
$t = \arccos \big( 1-\frac 1\alpha(1+\e^{- x})\big)$ in the
second.
We also rewrite $G$ by substituting $t =\arccos (1-x/\alpha)$
in the integral from \eqref{Def G}. Thus, with the abbreviations 
\[
 \beta &:=& 2\alpha-1\,\in\,\,]0,1] \\ 
 x_\alpha &:=& \, - \log \beta \,\in\,[0,\infty[ \\
 f(x) &:=&  
 \frac {\e^{- x}}{\pi\alpha} \sqrt{
 \frac {\beta +\e^{- x}}{1-\e^{- x} }} 
  +
 \frac {\e^{- x}}{\pi\alpha} \sqrt{
 \frac {\beta -\e^{-x}}{1+\e^{- x}}} 
 \1^{}_{]x_\alpha,\infty[}(x)  \\ 
g(x) &:=& \frac 1{\pi\alpha} \sqrt{\frac{2\alpha-x}x}\1^{}_{]0,2\alpha[}(x)
\]
for $x\in\,]0,\infty[$, we get
\la                         \label{F,G Laplace}
 F(\lambda,\alpha) 
\,=\, \int_0^\infty\exp\big(-\frac {\lambda x}\alpha \big)f(x)\d x,
&& 
 G(\lambda) \,=\, \int_0^\infty \exp\big(-\frac{\lambda x}\alpha\big)g(x)\d x 
\al
By letting $\lambda$ tend to 0 in  
\eqref{Def F alpha},  \eqref{Def G} and \eqref{F,G Laplace}, 
we conclude  that 
$\int_0^\infty f\d x = \int_0^\infty g\d x =1$.
Let us assume the following claim for the moment.

\medskip
{\sc Claim.} {\sl There is a $y_\alpha\in \,]0,\infty[$ with
$f<g$ on $]0,y_\alpha[$ and $f>g$ on $]y_\alpha,\infty[$. }

\medskip
Then, with $\phi(x):= \exp \big(- \lambda x/\alpha\big)$ being 
strictly decreasing in $x\in\,]0,\infty[$, we get 
$\int_0^\infty \phi f\d x < \int_0^\infty \phi g\d x$, and 
hence \eqref{crucial ineq}, by integrating  the inequality
$\big( \phi(x)- \phi(y_\alpha)\big) \big(f(x)-g(x) \big) <0$
over $x \in \,]0,\infty[\setminus\set{y_\alpha}$.  

Thus it only remains to prove the claim, which will be reduced
to two subclaims. Let us agree to use  interval notation
like $]a,b]:=\set{x\in\R : a<x \le b}$ also if $a=b$ and even if $a>b$. 

\medskip
{\sc Subclaim 1.} {\sl If $x_\alpha >0$, then 
there is a $u_\alpha\in \,]0,x_\alpha]$ with
$f<g$ on $]0,u_\alpha[$ and $f>g$ on $]u_\alpha,x_\alpha]$. }

\medskip
{\sc Subclaim 2.} {\sl There is a $v_\alpha\in [x_\alpha,\infty[$ with
$f<g$ on $[x_\alpha,v_\alpha[\setminus\set{0}$ and $f>g$ on 
$]v_\alpha,\infty[$. }

\medskip
{\sc Proof of the claim assuming the two subclaims.}
If $x_\alpha=0$, then we put $y_\alpha := v_\alpha$ and
observe that here $v_\alpha >0$, since both $f$ and $g$
integrate to $1$ over $]0,\infty[$. So let us assume 
now that $ x_\alpha>0$. 
For  $u_\alpha < x_\alpha$,  
Subclaim 1 yields $f(x_\alpha)>g(x_\alpha)$,
while for $v_\alpha> x_\alpha$,
Subclaim 2  yields $f(x_\alpha)<g(x_\alpha)$.  
Hence $u_\alpha = x_\alpha$
or $v_\alpha= x_\alpha$. We put 
 $y_\alpha := v_\alpha$ in the first case and
$y_\alpha := u_\alpha$ in the second. 
This proves the claim.

\medskip
{\sc Proof of Subclaim 1.}  
For  $x\in \,]\min(x_\alpha,2\alpha), x_\alpha]$,
we have $f(x)-g(x)=f(x)>0$, while for 
$x\in\,]0,\min(x_\alpha,2\alpha)]$, the difference 
$f(x)-g(x)$ has the same sign as 
\[
 \eta(x) &:=& \frac {\left( f^2(x) - g^2(x) \right) 
         (\pi\alpha)^2 x\e^{2x}(1-\e^{-x})}
 {\e^{2x} -\e^x - x } 
 \,=\,\frac {x(\e^{2x} -\e^x +\e^{-x} -1)}{\e^{2x} -\e^x - x } - 2\alpha
\] 
Differentiating $\eta$ according to the quotient rule yields a fraction 
with a positive denominator and with numerator 
\[
\e^{4x} -2\e^{3x} +(-2x^2+2x)\e^{2x} 
    +(x^2-4x+2)\e^{x} +2x-1+x^2\e^{-x} &=& 
   \sum_{k=2}^\infty \frac{(4x)^k}{k!} a_k \,>\,0  
\]
where  $a_2=a_3=0$,  while for   $k\ge 4$
\[
 a_k &:=& 1-2\left(\frac 34\right)^k -\frac {k(k-3)}{2^{k+1}}
      +\frac{k(k-1)(1+(-1)^k) -4k+2 }{4^k} \\
 &>& 1-2 \left(\frac 34\right)^k    
       -\frac{(k-1)(k-2)}{2^{k+1}} -\frac{k}{4^{k-1}} \,=:\,b_k 
 \,\ge\,  b_4 \,=\, \frac{15}{128} \,>\, 0
\]
by the monotonicity of the terms in the definition of $b_k$.
Thus $\eta$ is strictly increasing, and since   
$\lim_{x\rightarrow  0}\eta(x)=-2\alpha <0$,
Subclaim 1 follows with 
$u_\alpha:= \sup \{x\in\, ]0,\min(x_\alpha,2\alpha)]\,:\, \eta(x)<0 \}$. 

\medskip 
{\sc Proof of Subclaim 2.} 
For $x\in [x_\alpha, \infty[$
\[
 f_1(x) &:=& \left(\pi\alpha f(x)\sqrt{\e^{2 x}-1}\right)^2 
 \,=\,2  \left(\beta + \e^{-2 x} +
 \sqrt{(\beta^2 - \e^{-2 x})(1-\e^{-2 x})} \right)
\]
is increasing (possibly constant),
since for $x\in \,]x_\alpha, \infty[$
\[
 f_1'(x) &=& -4\e^{-2 x}
 \Big(1 - \frac {\frac 12(1+\beta^{2})-\e^{-2 x}}
    {\sqrt{(\beta^2 -\e^{-2 x})(1-\e^{-2 x})}}\Big) \,\ge \,0 
\]
by $\sqrt{ab} \le (a+b)/2$ with $a=\beta^2 -\e^{-2x}$
and $b=1-\e^{-2x}$. Further, 
$\lim_{x\rightarrow \infty}f_1(x) =4\beta \le4\alpha$.

For $x\in \,]0,2\alpha [$, let
\[
 g_1(x) &:=& \left(\pi\alpha g(x)\sqrt{\e^{2 x}-1}\right)^2
 \,=\, \frac{\e^{2 x}-1}{x}(2\alpha-x) 
\,=\, 4\alpha + \sum_{k=1}^\infty \frac{(2x)^k}{k!}
        \Big(\frac{4\alpha}{k+1} -1\Big) 
\]
Then $g_1$ is strictly decreasing on the interval
$\set{x\in\,]0,2\alpha[\,:\, g_1(x) <4\alpha}$,  
since $\lim^{}_{x\downarrow 0} g_1(x)=4\alpha$
 and since $g_1'(x) =: \sum_{k=0}^\infty c_k x^k$ is, for some $\delta>0$,
strictly positive on $]0,\delta[$ and 
strictly negative for $x>\delta$, 
as the coefficients $c_k= (2^{k+1}/k!)\big(4\alpha/(k+2) -1\big)$
change sign exactly once, from plus to minus. See
 P\'olya \& Szeg\"o (1971, p.~43, Aufgaben V.38 and V.40) for
this last argument.

Hence there is at most one $x$ 
in the possibly empty interval  $[x_\alpha,2\alpha[$  with $f_1(x)=g_1(x)$.
We take this $x$, if it exists, as $v_\alpha$, and otherwise put 
$v_\alpha :=x_\alpha$. 
The subclaim follows: For 
$x\in [x_\alpha, 2\alpha[\setminus\set{0}$, 
the difference $f(x) -g(x)$ has the same  
sign as $f_1(x)-g_1(x)$, while for 
$x\ge 2\alpha$, we have $f(x)>0=g(x)$.  \Halmos  

\subsection{Proof of Theorem \ref{Kanter Lemma 4.3}} 
\label{Proof of  Kanter's Lemma 4.3}
By Remark \ref{Remarks in its own right}\ref{Where max?},
it suffices to prove the claim for $k=0$. Fourier inversion
yields 
\la
 \STPC_p(\set{0,1})&=&\frac 1\pii \int_0^\pi
 \Big(\prod_{j=1}^n \big(1-p_j(1-\cos t)\big) \Big)
   \big( 1+ \cos t)\d t \label{Qp01 via Fourier} 
\al
for $p\in\cP$.
Let us fix $n\in\N$ and $\lambda \in \,]0,n]$ for the rest of 
this proof. By Hoeffding (1956, Corollary 2.1), the probability
$ \STPC_p(\set{0,1})$, being a permutation invariant and multiaffine 
function of $p\in [0,1]^n$, attains its maximum subject to
the constraint $|p|= \lambda$ at some
point $(1,\ldots,1,\alpha,\ldots,\alpha,0,\ldots,0)\in [0,1]^n$
having  at most three different coordinate values,
with at most one of them distinct from $0$ and $1$. Thus,  
for every $p\in[0,1]^n$ with $|p| =\lambda$,
there exist $\ell,m \in\N_0$ and $\alpha \in\,]0,1[$
with $\ell + m\le n$, $\ell+m \alpha = \lambda$, and 
\newlength{\formulawidth}
\settowidth{
      \formulawidth}{$\displaystyle \max
\big(F(\lambda,\alpha),F(\lambda,1) \big)
   \qquad\qquad \qquad\qquad$}
\newcommand{\formulabox}[1]{
\makebox[\formulawidth][l]{$\displaystyle #1$}}
\la
 \nonumber
 \STPC_p(\set{0,1})&\le &\frac 1\pii \int_0^\pi
      \big( \cos t \big)^\ell\big( 1-\alpha(1-\cos t)\big)^m(1+\cos t)\d t
              \qquad \text{[using \eqref{Qp01 via Fourier}]}\\
 \nonumber
 &\le&\frac 1\pii \int_0^\pi
      \big| \cos t \big|^\ell
 \big| 1-\alpha(1-\cos t)\big|^{\frac{\lambda -\ell}\alpha }(1+\cos t)\d t\\
 &\le & \max \big(F(\lambda,\alpha),F(\lambda,1) \big)
   \qquad\qquad \qquad\qquad  \text{[see below and \eqref{Def F alpha}]}
  \label{convexity argument}\\
 \nonumber
 &<&\formulabox{G(\lambda)}
   \text{[by Lemma \ref{special analytic inequality}]}\\
 \nonumber
 &=&  \formulabox{\SymPois_\lambda(\set{0,1})}
   \text{[by \eqref{Glambda vs Slambda01}
]}
\al
Here \eqref{convexity argument} follows by regarding the integrand, 
and hence the integral, 
in the preceding line as a convex function of the parameter 
$\ell\in[0,\lambda]$, so that the integral becomes maximal for 
$\ell= 0$ or $\ell=\lambda$.\Halmos

\subsection{Proof of Lemma \ref{Lemma Simple bounds for G}} 
The Laplace transform representation  \eqref{F,G Laplace} with $\alpha=1$,
\la                 \label{G as Laplace}
  G(\lambda) \,=\, \int_0^\infty \e^{-\lambda x} g(x)\d x, &&
  g(x) \,=\,  \frac 1\pii\sqrt{\frac {2-x}x }
 \1^{}_{]0,2[}(x)   
\al
yields the complete monotonicity of $G$, 
compare Berg {\it et al.} (1984, p.~135).

It is easily checked that there 
is no need to qualify the big O claims in the lemma
or this proof by ``$\lambda \rightarrow \infty$'' or similarly.
Standard Bessel function asymptotics, see Olver (1997, pp. 251, 238),
yield for every $k\in \Z$ 
\la                   \label{BesselI for large argument}
 \SymPois_\lambda(\set{k}) &=& \mathrm{e}^{-\lambda}\BesselI_k(\lambda) 
 \,=\, \frac 1{\sqrt{2\pi \lambda}}
  \Big(1 -\frac{4k^2 - 1}{8\lambda}  
   + \frac{(4k^2 - 1)(4k^2 - 9)}{128\lambda^2}
   + O(\frac 1{\lambda^3}) \Big)
\al
and hence \eqref{G asymp for lambda large}.
Also, for fixed $k\in\N_0$,
a multiplication of power series yields 
\[
  \SymPois_\lambda(\set{k}) &=& 
  \Big(\sum_{j=0}^\infty\frac{(-\lambda)^j}{j!}\Big) 
  \Big(\sum_{j=0}^\infty\frac
  {(\lambda /2)^{2j+k}}{j!(j+k)!}\Big)
  \,=\, \frac{\lambda^{k}}{2^{k}k!}  
       -   \frac{\lambda^{k+1}}{2^{k}k!} 
   + \frac{(2k+3)\lambda^{k+2}}{2^{k+2}(k+1)!}
   + O(\lambda^{k+3})
\]
and hence \eqref{G asymp for lambda small}.
The big O claims in \eqref{simple bound, 2}
and \eqref{simple bound, 3} are obvious from the
binomial series. 

For $a,b\in [0,\infty[$, we have 
\la        \label{roots as Laplace transforms}
 \sqrt{\frac a{b+\lambda}} 
  &=& \int_0^\infty \mathrm{e}^{-\lambda x}\sqrt{\frac a{\pii x}}
     \mathrm{e}^{-bx}\d x 
\al
Taking here $a=2/\pi$ and $b=1/4$ and using 
$\mathrm{e}^{-x/4} =\sqrt{\mathrm{e}^{-x/2}} 
> \sqrt{1 - \frac x2}\1^{}_{]0,2[}(x)$
and hence 
$\frac 1\pii \sqrt{\frac 2x}\mathrm{e}^{-x/4} >g(x)$,
compare \eqref{G as Laplace},
yields  the inequality in \eqref{simple bound, 2}. 
Inequality \eqref{simple bound, 1} follows trivially.
Taking $a=b=1$ in \eqref{roots as Laplace transforms} yields
\[
  H(\lambda) &:=& \frac 1{\sqrt{1+\lambda}} 
    \,=\, \int_0^\infty \mathrm{e}^{-\lambda x} h(x)\d x
\] 
with $h(x):= \mathrm{e}^{-x}/\sqrt{\pii x}$.
Now  
\la       \label{masses and expectations equal} 
   \int_0^\infty h(x)\d x = 1 = \int_0^\infty g(x)\d x
 &\text{and}& \int_0^\infty xh(x)\d x = \frac 12  = \int_0^\infty xg(x)\d x
\al
On $]0,2[$, the derivative 
$\big( \log(g(x)/h(x))\big)' = \big(\frac 32 -x \big)/(2-x)$
changes sign exactly once, from plus to minus.
Further, $\lim_{x\downarrow 0} g(x)/h(x) = \sqrt{2/\pii} <1$, 
$g/h$ is continuous on $]0,\infty[$, and $g(x)=0<h(x)$ for $x\in [2,\infty[$.
Recalling the first half of  \eqref{masses and expectations equal},
we deduce the existence of $x_1,x_2\in\,]0,2[$ with $x_1<x_2$,
$g<h$ on $]0,x_1[ \,\cup\, ]x_2,\infty[$, 
and $g > h$ on $]x_1,x_2[$. Thus, since $\phi(x):= \mathrm{e}^{-\lambda x}$
is strictly convex, we have
\[
 \Big(\phi(x)- \big(\phi(x_1)\frac{x_2-x}{x_2-x_1} 
  + \phi(x_2) \frac {x-x_1}{x_2-x_1}   \big) \Big)\big(h(x)-g(x)\big)
 &>& 0
\] 
for $x\in \, ]0,\infty[\setminus \set{x_1,x_2}$,
and hence obtain $G(\lambda) < H(\lambda)$ by integration,
using \eqref{masses and expectations equal}. \Halmos

\subsection{Proof of Corollary \ref{Corollary improving Roos 2005}}
Clearly $\otherp_j\geq 1-\Conc(X_j,t)$ for every $j$.
Hence 
\[
 \Conc(\sum_{j=1}^n X_j, t) \,\le\, \BerC_\otherp \psi
  \,=\,  \STP_\otherp(\set{0,1}) \,\le\, G(|\otherp|) 
  \,\le\, G\Big(\sum_{j=1}^n \big(1-\Conc(X_j,t)\big)\Big)
\]  
by Le~Cam (1986, p.~411, Proof of Theorem 2),
\eqref{ConcQp as Bernoulli conv expectation},
\eqref{Kanter Lemma 4.3 inequality}, \eqref{Glambda vs Slambda01},
and the antitonicity of $G$ from Lemma~\ref{Lemma Simple bounds for G}.
Applying \eqref{simple bound, 2} and \eqref{simple bound, 3}
yields the claimed inequalities.\Halmos

\subsection{Proof of Corollary \ref{Corollary improving Barbour Xia 1999}}
By using the Mineka coupling, Barbour \& Xia~(1999, Proposition~4.6)
proved that
\[
d_{\mathrm{TV}}  (\sum_{j=1}^n X_j,\, 1+\sum_{j=1}^n X_j)
 &\le& \P(\sum_{j=1}^nY_j\in\{0,1\})
\]
where the $Y_1,\dots,Y_n$ are independent and symmetric 
as in \eqref{Distr of Yj}, but now with 
$ p_j  = 1-d_{\mathrm{TV}}(X_j,1+X_j)$ for $j=1,\dots,n$.
Thus  Theorem~\ref{Theorem (Kanter's Lemma 4.3)}
and (\ref{simple bound, 2}) yield our claim.
\Halmos

\mbox{}

\medskip {\bf Acknowledgements.}
We  thank Uwe R\"osler for  discussions prompting 
Remark~\ref{Remarks in its own right}\ref{Remark on Marshall and Olkin}.  
Bero Roos thanks Professor Wilfried Grecksch of 
Martin--Luther-University  Halle--Witten\-berg
for an invitation to  the most hospitable 
Institute of Optimization and Stochastics 
for the  Winter term 2005/06, during which this work was finished.

\bigskip
{\Large\bf References } 

\medskip
{\footnotesize
\begin{description}

\item[\sc Barbour, A.D., Holst, L. \&  Janson, S. (1992).]
 {\it Poisson Approximation.} Clarendon Press, Oxford.

\item[\sc Barbour, A.D. \& Xia, A. (1999).]  Poisson perturbations.
{\it ESAIM Probab.\ Statist.} {\bf 3}, 131--150.
 
\item[\sc Berg, C., Christensen, J.P.R. \& Ressel, P. (1984).]
{\it Harmonic Analysis on Semigroups.} Springer, New York.

\item[\sc Bickel, P.J.  \& van Zwet, W.R. (1980).]
On a theorem of Hoeffding. 
In: {\it Asymptotic Theory of Statistical Tests and Estimation:
in Honor of Wassily Hoeffding, I.M. Chakravarti (ed.), 
Academic Press, New York, pages 307--324.}

\item[\sc Boland, P.J. \& Proschan, F. (1983).]
The reliability of $k$ out of $n$ systems.
{\it Ann.\ Probab.} {\bf 11}, 760-764.

\item[\sc Bondar, J.V. (1994).]
Comments and complements to 
{\it Inequalities: Theory of Majorization and Its Applications}
by Albert W.\ Marshall and Ingram Olkin.
{\it Linear Algebra and its Applications} {\bf 199}, 115--130.

\item[\sc Bretagnolle, J. (2004).] 
Sur l'in\'egalit\'e de concentration de Doeblin-L\'evy, Rogozin-Kesten. 
In: {\it 
Parametric and semiparametric models with applications to reliability,
survival analysis, and quality of life, 533--551,
Stat. Ind. Technol., Birkh\"auser, Boston.}

\item[\sc \v{C}ekanavi\v{c}ius, V.   \& Roos, B. (2006).]
An expansion in the exponent for compound binomial approximations.
{\it Lith.\ Math.\ J.} {\bf 46}, 54--91.

\item[\sc Dharmadhikari, S. \& Joag-Dev, K. (1988).]
{\it Unimodality, Convexity, and Applications.} Academic Press, Boston.

\item[\sc Gleser, L. (1975).] On the distribution of the number of successes
in independent trials. {\it Ann.\ Probab.} {\bf 3}, 182-188.
\item[\sc Hoeffding, W. (1956).] 
On the distribution of the number of successes in independent trials.
{\it Ann.\ Math.\ Statist.} {\bf 27}, 713--721.
Also in: {\it The Collected Works of Wassily Hoeffding},
 N.I.\ Fisher \& P.K.\ Sen (eds.), Springer, New York, 1994.

\item[\sc Kanter, M. (1976).] Probability inequalities for convex sets
and multidimensional concentration functions.
{\it J. Multivariate Anal.} {\bf 6}, 222--236.  

\item[\sc Karlin, S. \& Novikoff, A. (1963).]
Generalized convex inequalities.
{\it Pacific J.\ Math.} {\bf 13}, 1251--1279.  

\item[\sc Le~Cam, L. (1986).]
  {\it Asymptotic Methods in Statistical Decision Theory.}
  Springer, New York.

\item[\sc Marshall, A.W. \& Olkin, I. (1979).]
{\it Inequalities: Theory of Majorization and Its Applications.}
Academic Press, New York. 

\item[\sc Mattner, L. (2006).]
Lower bounds for tails of sums of independent symmetric random variables. 
Preprint, \verb§http://arxiv.org/abs/math.PR/0609200§.

\item[\sc Merkle, M. \& Petrovi\'c, L. (1997).] 
Inequalities for  sums independent geometrical random variables.
{\it Aequ.\ Math.} {\bf 54}, 173--180. 

\item[\sc  Olver, F.W.J. (1997).] 
{\it Asymptotics and Special Functions.}
A K Peters, Wellesley, Massachusetts.

\item[\sc Petrov, V.V. (1995).] 
{\it Limit Theorems of Probability Theory. 
 Sequences of Independent Random Variables.}
 Clarendon Press, Oxford.    

\item[\sc P\'olya, G. \& Szeg\"o, G. (1971).]
{\it Aufgaben und Lehrs\"atze aus der Analysis II.}
4.\ Auflage, Springer, Berlin.

\item[\sc Rogozin, B.A. (1993).]
Inequalities for concentration of a decomposition.
{\it Theory Probab.\ Appl.} {\bf 38}, 556--562.

\item[\sc Roos, B. (2005).] 
  On Hipp's compound Poisson approximations 
  via concentration functions. 
  {\it Bernoulli} {\bf 11}, 533--557.

\end{description}
}
{\footnotesize
\begin{multicols}{2}
\noindent
{\sc Universit\"at zu L\"ubeck\\ 
     Institut f\"ur Mathematik\\
     Wallstr.\ 40\\
     D-23560 L\"ubeck\\
     Germany\\
     Email: \verb§mattner@math.uni-luebeck.de§ \\
     Universit\"at Hamburg \\
     Department Mathematik, SPST  \\
     Bundesstr.\ 55  \\
     D-20146 Hamburg  \\
     Germany \\
     Email:  \verb§roos@math.uni-hamburg.de§
}  
\end{multicols}
}
\end{document}